\documentclass[12pt]{amsart}
\usepackage{latexsym,amsfonts,amsthm,amsmath,amscd,amssymb}

\usepackage[all]{xy}
\usepackage{array}

\newtheorem{theorem}{Theorem}
\newtheorem{proposition}[theorem]{Proposition}
\newtheorem{lemma}[theorem]{Lemma}
\newtheorem*{conjecture}{Generalized Chen's Conjucture}

\begin{document}

\title{New examples of biharmonic submanifolds in $CP^n$ and
$S^{2n+1}$}

\author{wei zhang}

\thanks{Mathematics Classification Primary(2000): 58E20.}

\maketitle

\begin{abstract}

We construct biharmonic real hypersurfaces and Lagrangian
submanifolds of Clifford torus type in $CP^n$ via the Hopf
fibration; and get new examples of biharmonic submanifolds in
$S^{2n+1}$ as byproducts .

\end{abstract}

\section{Introduction}
Let $\phi: (M,g)\rightarrow(N,h)$ be a smooth map between two
Riemannian manifolds and define its first tension field as
$\tau(\phi)=trace \nabla d\phi$. $\phi$ is called harmonic if it is
a critical point of the energy functional:
$$E(\phi)=\frac{1}{2}\int_M |d\phi|^2 v_g$$
This amounts to $\tau(\phi)=0$. In the isometric immersion case,
harmonic is equivalent to minimal.

Considering the bienergy $E_2(\phi)=\frac{1}{2}\int_M |\tau(\phi)|^2
v_g$, its critical points are define as biharmonic maps (\cite
{BB}). The associated Euler-Lagrange equation is given by the
vanishing of the bitension field:
$$\tau_2(\phi)=-\Delta^\phi \tau(\phi)-trace
R^N(d\phi,\tau(\phi))d\phi$$

Obviously, any harmonic map is biharmonic. We call the non-harmonic
one proper biharmonic. A submanifold is called biharmonic if the
inclusion map is biharmonic.

Concerning the proper biharmonic map, there are several
non-existence results for the non-positive sectional curvature
codomains(\cite{J1,BB}), for instance:

\begin{theorem}[\cite{J1,J2}]
Let $\phi: (M,g)\rightarrow (N,h)$ be a smooth map. If M is compact,
orientable and $Riem^N\leq0$, then $\phi$ is biharmonic if and only
if it is harmonic.
\end{theorem}

Suggested by this kind of results:
\begin{conjecture}
Biharmonic submanifolds of a manifold N with $Riem^N\leq0$ are
minimal.
\end{conjecture}

So it is sensible to focus on the proper biharmonic submanifolds in
sphere or other non-negatively curved space. Although some partial
classification results had been obtained (\cite{BMO}), the known
examples are still relatively rare. In this article, we will
construct series new examples in $S^{2n+1}$ and $CP^n$.

\begin{theorem}
Denote $M^C_{p,q}(r,s)$ be the image of the generalized Clifford
torus $M_{p,q}(r,s)=S^{2p+1}(r)\times S^{2q+1}(s)$ in $S^{2n+1}$ via
the Hopf map, where p+q=n-1. Then it is biharmonic in $CP^n$ if and
only if:
\begin{equation}
(\frac{s}{r})^2(2p+1)+(\frac{r}{s})^2(2q+1)=2(n+2)
\end{equation}

\end{theorem}
And

\begin{theorem}
Identify $R^{2n+2}$ with $C^{n+1}$. Define the Clifford torus
$T^{n+1}=\{|z_i|=a_i|\sum^{n+1}_{i=1}a_i^2=1\}$ in $S^{2n+1}$. Then
it is biharmonic if and only if:

\begin{equation}
a_id-\frac{1}{a_i^3}=2(n+1)((n+1)a_i-\frac{1}{a_i})
\end{equation}
where $d=\sum \frac{1}{a_i^2}$, $i=1,\ldots, n+1$.

\end{theorem}

Followed by:

\begin{theorem}
Let $T_C^n$ be the quotient of $T^{n+1}$ by the $S^1$ action. Then
it is biharmonic in $CP^n$ if and only if:

\begin{equation}
a_id-\frac{1}{a_i^3}=2(n+3)((n+1)a_i-\frac{1}{a_i})
\end{equation}

\end{theorem}

The explicit expressions of $r$,$s$ and $a_i$ will be solved out in
section 3.

Among these torus, we can pick out plenty proper ones. As the
construction is very routinely, we would emphasis on the method
rather than the concrete examples.

\section{Preliminary}
\subsection{Hopf fibration}

There are already many known examples of biharmonic submanifolds in
sphere. To find examples in $CP^n$, the Hopf fibration: $\pi:
S^{2n+1}\xrightarrow{S^1} CP^n$ is a natural candidate, where $\pi$
is a Riemann submersion with totally geodesic fibres $S^1$ and
$CP^n$ has constant holomorphic sectional curvature 4.

This submersion has lots of good properties, such as:

\begin{lemma}[\cite{WL}]
If $\tilde{M}$, M are submanifolds respect this Riemann submersion,
i.e. the following diagram commutes:

$$
\xymatrix{
  \tilde{M} \ar[d]_{\pi} \ar[r]
                & S^{2n+1} \ar[d]^{\pi}  \\
  M   \ar[r]
                & CP^n             }
$$
then $\tilde{M}$ minimal(totally geodesic) is equivalent to M
minimal(totally geodesic). More precisely, the mean curvature
$\tilde{H}$ of $\tilde{M}$ is the horizontal lift of H of M.
\end{lemma}

In this paper, to simplify the notations, we take the convention
that $H$ is non-normalized, i.e $H=\tau(i)$.

Unfortunately, this submersion do not preserve the biharmonicity.
See \cite {LO}, consider the Hopf map $S^3\rightarrow
S^2(\frac{1}{2})$. Lifting the biharmonic submanifold
$S^1(\frac{\sqrt{2}}{4})$, we get the Clifford tori
$S^1(\frac{\sqrt{2+\sqrt{2}}}{2})\times
S^1(\frac{\sqrt{2-\sqrt{2}}}{2})$, which is not biharmonic in $S^3$.

But it provides us the stereotype of constructing biharmonic
submanifolds in $CP^n$: If we modify the radius of the standard
Clifford torus properly, its image would be biharmonic in $CP^n$.

\subsection{The biharmonic equations in $S^n$ and $CP^n$}
When the submanifold lies in $S^n$, it is convenient to split the
bitension field in its normal and tangent components.

\begin{theorem}[\cite {CMO}]
$M^m$ is a submanifold of $S^n$, then it is biharmonic if and only
if
\begin{equation}{}
\begin{cases}
-\Delta^{\bot}H-traceB(\cdot, A_H \cdot)+mH=0 & {} \\
2traceA_{\nabla^{\bot}_{(\cdot)}H}(\cdot)+\frac{1}{2}grad(|H|^2)=0
\end{cases}
\end{equation}
Moreover, if we assume M has parallel mean curvature, the equation
turns to:

\begin{equation}
\sum_{i,j}B_{ij}<B_{ij},H>=mH
\end{equation}
\end{theorem}
Where A denotes the Weingarten operator, B the second fundamental
form, $\nabla^\bot$ and $\Delta^\bot$ the connection and the
Laplacian in the normal bundle.

When the ambient space is $CP^n$, considering the real hypersurfaces
or the Lagrangian submanifolds, the biharmonic equation has similar
form:

\begin{proposition}
$M^{2n-1}$ is a real hypersurface in $CP^n$, then it is biharmonic
if and only if
\begin{equation}{}
\begin{cases}
-\Delta^{\bot}H-traceB(\cdot, A_H \cdot)+2(n+1)H=0 & {} \\
2traceA_{\nabla^{\bot}_{(\cdot)}H}(\cdot)+\frac{1}{2}grad(|H|^2)=0
\end{cases}
\end{equation}
If M has parallel mean curvature in addition, the equation becomes:

\begin{equation}
||B||^2=2(n+1)
\end{equation}

\end{proposition}

And

\begin{proposition}
$M^n$ is a Lagrangian submanifold of $CP^n$, then it is biharmonic
if and only if
\begin{equation}{}
\begin{cases}
-\Delta^{\bot}H-traceB(\cdot, A_H \cdot)+(n+3)H=0 & {} \\
2traceA_{\nabla^{\bot}_{(\cdot)}H}(\cdot)+\frac{1}{2}grad(|H|^2)=0
\end{cases}
\end{equation}
When M with parallel mean curvature, it is simplified to:

\begin{equation}
\sum_{i,j}B_{ij}<B_{ij},H>=(n+3)H
\end{equation}

\end{proposition}

\noindent \emph{Proof of Proposition 7:} Denote the canonical
inclusion map as $i$, then $\tau(i)=H$.

 Since $traceR^{CP^n}(di,\tau(i))di=-2(n+1)\tau(i)=-2(n+1)H$,
 $i$ biharmonic $\Leftrightarrow$
 $$\tau_2(i)=trace\nabla dH+2(n+1)H=0$$
 Choose normal frame $\{e_{\alpha}\}$ on M, where $1\leq \alpha \leq 2n-1$.
 \begin{equation*}
 \begin{split}
 &trace \nabla dH=
 \nabla_{e_\alpha}^{CP^n}\nabla_{e_\alpha}^{CP^n}H=\nabla_{e_\alpha}^{CP^n}(\nabla_{e_\alpha}^\bot
 H-A_H(e_\alpha))\\
 =&\nabla_{e_\alpha}^\bot \nabla_{e_\alpha}^\bot
 H-A_{\nabla_{e_\alpha}^{\bot}H}(e_\alpha)-\nabla_{e_\alpha}A_H(e_\alpha)-B(e_\alpha,
 A_H (e_\alpha))\\
 =&-\Delta^{\bot}H-traceB(\cdot, A_H
 \cdot)-(A_{\nabla_{e_\alpha}^{\bot}H}(e_\alpha)+\nabla_{e_\alpha}A_H(e_\alpha))
 \end{split}
 \end{equation*}

 Rewrite $A_H$ and use Codazzi equation:
 \begin{equation*}
 \begin{split}
&\nabla_{e_\alpha}A_H(e_\alpha)=(\nabla_{e_\alpha}<B(e_\beta,
e_\alpha),H>)e_\beta\\
=&<B(e_\beta,
e_\alpha),\nabla_{e_\alpha}^{\bot}H>e_\beta+<(\nabla_{e_\alpha}^{\bot}B)(e_\beta,
e_\alpha),H>)e_\beta\\
=&A_{\nabla_{e_\alpha}^{\bot}H}(e_\alpha)+<(\nabla_{e_\beta}^{\bot}B)(e_\alpha,
e_\alpha),H>)e_\beta+<R^{CP^n}(e_\alpha,e_\beta)e_\alpha,H>e_\beta\\
=&A_{\nabla_{e_\alpha}^{\bot}H}(e_\alpha)+\frac{1}{2}grad(|H|^2)
 \end{split}
 \end{equation*}
where
$<R^{CP^n}(e_\alpha,e_\beta)e_\alpha,H>=Ric^{CP^n}(e_\beta,H)=0$ for
any $\beta$ by $CP^n$ is Einstein.

Replacing $trace \nabla dH$ in the identity and arranging the terms
in tangent and normal components, we get what desired.

 \qed

Proof of Proposition 8 is similar.

Our general modus operandi is adjusting the radius of the Clifford
torus until their second fundamental forms satisfying the biharmonic
equations, and the fact that the Clifford type torus we concern
about all have parallel mean curvature simplifies the computation.

\section{The examples}
\subsection{Generalized circles $M^C_{p,q}(r,s)$ in $CP^n$}
In \cite{L}, Lawson introduced the concept of generalized equator
$M^C_{p,q}$ which is mininal in $CP^n$. Thus we call
$M^C_{p,q}(r,s)$ generalized circle which is not always minimal.

For the generalized Clifford torus in $S^{2n+1}$ have constant mean
curvature, so does $M^C_{p,q}(r,s)$.

By (7), $M^C_{p,q}(r,s)$ is biharmonic if $||B||^2=2(n+1)$.

While form \cite{L}, we know $||B||^2=||\tilde{B}||^2-2$, where
$\tilde{B}$ is the second fundamental form of the generalized
Clifford torus in $S^{2n+1}$. Combining the fact $||\tilde{B}||^2=
(\frac{s}{r})^2(2p+1)+(\frac{r}{s})^2(2q+1)$ leads to (1). It is
easy to check the solution is never minimal.

This is a table of biharmonic real hypersurfaces in $CP^n$(We only
list the case n=5).

\begin{tabular}{|c|c|c|c|c|}
  \hline
  n & p & q & r                            & s    \\
   \hline
  5 & 0 & 4 & $\sqrt{(3*(5 - \sqrt3))/22}$ & $\sqrt{(7 +(3*\sqrt3))/22}$ \\
   \hline
  5 & 0 & 4 & $\sqrt{(3*(5 + \sqrt3))/22}$ & $\sqrt{(7 -(3*\sqrt3))/22}$ \\
   \hline
  5 & 1 & 3 & $\sqrt{(13 - \sqrt{15})/22}$   & $\sqrt{(9 +
  \sqrt{15})/22}$\\
   \hline
  5 & 1 & 3 & $\sqrt{(13 + \sqrt{15})/22}$   & $\sqrt{(9 - \sqrt{15})/22}$\\
   \hline
  5 & 2 & 2 & $\sqrt{(11 + \sqrt{11})/22}$   & $\sqrt{(11 - \sqrt{11})/22}$\\

  \hline
\end{tabular}

\subsection{Clifford torus $T^{n+1}$ in $S^{2n+1}$}

Denote the position vector by $x=a_ix_i$, where $x_i$ is the unit
vector. Choosing normal frame $\{e_i\}$ of the tangent space s.t
$Je_i=x_i$ where $J$ is the complex structure in $C^{n+1}$. Direct
computation show that $B_{ij}=\delta_{ij}(-\frac{x_i}{a_i}+x)$ and
$H=\sum(a_i(n+1)-\frac{1}{a_i})x_i$.

Thus
$$\sum B_{ij}<B_{ij},H>=(\sum(a_id-\frac{1}{a_i^3})x_i)-(n+1)H$$
where $d=\sum_{i=1}^{n+1}\frac{1}{a_i^2}$.

Furthermore, $T^{n+1}$ has parallel mean curvature. If it is
biharmonic, by (5)
$$(\sum(a_id-\frac{1}{a_i^3})x_i)-(n+1)H=(n+1)H$$

Compare each components of H, we have (2).

This equation system is some how funny. Denote $b_i=a_i^2$, rewrite
it as:

\begin{equation}
 (2(n+1)^2-d)b_i^2-2(n+1)b_i+1=0
\end{equation}
where $d=\sum_{i=1}^{n+1}\frac{1}{b_i}$.

There are (n+2) equations and (n+1) variables, while the equations
are not independent. It is always solvable. To see this, view d as a
known number. Let $r$, $s$ be the two different roots of the
quadratic (10), and $p$, $q$ be their multiplicities in $b_i$
respectively. Then $pr+qs=1$, and

\begin{equation*}
d=\frac{p}{r}+\frac{q}{s}=\frac{p}{r}+\frac{p}{s}+\frac{q-p}{s}=p\frac{r+s}{rs}+\frac{q-p}{\frac{1-p(r+s)}{q-p}}
\end{equation*}

Let $t=2(n+1)^2-d$, use Vieta's Theorem:

$$(2q(n+1)-t)(t-2p(n+1))=(q-p)^2t$$
i.e.

$$t^2-(2(n+1)^2-(p-q)^2)t+4pq(n+1)^2=0$$

This equation's discriminant is:
\begin{equation*}
\begin{split}
\Delta &=(p-q)^4+4(n+1)^4-(4(p-q)^2(n+1)^2+16pq(n+1)^2)\\
&=(p-q)^4+4(n+1)^4-4(p+q)^2(n+1)^2=(p-q)^4
\end{split}
\end{equation*}

The root $t=(n+1)^2$ leads to minimal torus, so we pick the one
$t=(n+1)^2-(p-q)^2$. The original quadratic becomes:

$$((n+1)^2-(p-q)^2)b_i^2-2(n+1)b_i+1=0$$
Its two roots are $\frac{1}{(n+1)\pm(p-q)}$ with multiplicities p
and q.

Thus the torus is with elegant form
$T^{n+1}_{p,q}(\frac{1}{\sqrt{2p}},\frac{1}{\sqrt{2q}})$. It can not
be minimal unless p=q.

\subsection{Clifford torus $T_C^n$ in $CP^n$}

$T_C^n$ is the quotient of $T^{n+1}$ by the $S^1$ action. It is
Lagrangian in $CP^n$, and has parallel mean curvature fields.

Choose orthonormal frame $\{f_i\}$,$i=1,\ldots,n$. To compute
$\sum_{i,j}B_{i,j}<B_{i,j},H>$, as the Hopf map preserve the
horizontal part of second fundamental(\cite{WL}) and the mean
curvature, we lift it to $S^{2n+1}$. Without ambiguity, denote
$\{f_i,\nu\}$ the orthonormal basis of $T^{n+1}$, where $\nu=Jx$,
$x$ the position vector, is the vertical vector fields. Notice that
$\sum_{\alpha,\beta}\tilde{B}_{\alpha,\beta}<\tilde{B}_{\alpha,\beta},\tilde{H}>$
is independent of basis, thus:
\begin{equation*}
\begin{split}
&\sum_{\alpha,\beta}\tilde{B}_{\alpha,\beta}<\tilde{B}_{\alpha,\beta},\tilde{H}>-2\tilde{B}_{i,\nu}<\tilde{B}_{i,\nu},\tilde{H}>-\tilde{B}_{\nu,\nu}<\tilde{B}_{\nu,\nu},\tilde{H}>\\
=&\widetilde{\sum_{i,j}B_{i,j}<B_{i,j},H>}
\end{split}
\end{equation*}
It is not hard to show:

$$\tilde{B}_{i,\nu}<\tilde{B}_{i,\nu},\tilde{H}>=-Jf_i<-Jf_i,\tilde{H}>=\tilde{H}$$
and $\tilde{B}_{\nu,\nu}=0$.

Follow the biharmonic equation (9):

$$\sum_{\alpha,\beta}\tilde{B}_{\alpha,\beta}<\tilde{B}_{\alpha,\beta},\tilde{H}>=(n+5)\tilde{H}$$

This is nothing but (3). We solve it in the same way.

\begin{equation*}
d=\frac{p}{r}+\frac{p}{s}+\frac{q-p}{s}=p\frac{r+s}{rs}+\frac{q-p}{\frac{1-p(r+s)}{q-p}}
\end{equation*}

Let $t=2(n+1)(n+3)-d$, we have:
\begin{equation}
t^2-(2(n+1)(n+3)-(p-q)^2)t+4pq(n+3)^2=0
\end{equation}

Its discriminant is:

$$\Delta=(p-q)^4+8(n+3)(p-q)^2$$

Both roots of t make the original quadratic about r and s solvable.
Same as the generalized circles $M^C_{p,q}(r,s)$ case, the
biharmonic torus ${T^n_C}_{p,q}$ could never be minimal.

The software "Mathematica" is helpful to get the explicit
expressions. The following is the table of ${T^n_C}_{p,q}(r,s)$ when
n=4:

\begin{tabular}{|c|c|c|c|c|}
  \hline
  n & p & q & r                            & s    \\
   \hline
  4 & 1 & 4 & $\sqrt{(11 - \sqrt{65})/7}/2$ & $\sqrt{(17 + \sqrt{65})/7}/4$ \\
   \hline
  4 & 1 & 4 & $\sqrt{(11 + \sqrt{65})/7}/2$ & $\sqrt{(17 - \sqrt{65})/7}/4$ \\
   \hline
  4 & 2 & 3 & $\sqrt{(13 - \sqrt{57})/14}/2$& $\sqrt{(15 + \sqrt{57})/21}/2$\\
  \hline
  4 & 2 & 3 & $\sqrt{((13+\sqrt{57})/14)}/2$&$\sqrt{(15 -
  \sqrt{57})/21}/2$\\
  \hline
\end{tabular}
\\

We add a remark on the biharmonic stability of such type torus. A
biharmonic submanifold is called stable if for any variation of the
inclusion map $\{i_t\}$, its second derivative of $E_2$ is
nonnegative. Using $H=\tau(i)$ as the variation vector fields, by
the second variation formula in \cite{J2}:

\begin{equation*}
\begin{split}
\frac{1}{2}\frac{d^2}{dt^2}E_2(i_t)|_{t=0}=&4\int_T
<R^{CP^n}(f_k,H)\nabla_{f_k}^{CP^n}H,H>d v_g\\
=&4\int_T(-||H||^4+3<f_k,JH><\nabla_{f_k}^{CP^n}H,JH>)d v_g
\end{split}
\end{equation*}

Notice that $JH=<JH,f_l>f_l$ and
$$<\nabla_{f_k}^{CP^n}H,f_l>=-<\nabla_{f_k}^{CP^n}f_l,H>=-<B_{kl},H>$$
the above expression becomes:

\begin{equation}
-4\int_T(||H||^4+3<Jf_k,H><B_{kl},H><Jf_l,H>)dv_g
\end{equation}

Still lift them to $S^{2n+1}\subset C^{n+1}$ for computation. Using
$\{f_i,\nu\}$ and $\{e_\alpha\}$ as frames alternatively, as
$$<Jf_k,H><B_{kl},H><Jf_l,H>=<\tilde{Jf_k},\tilde{H}><\tilde{B}_{kl},\tilde{H}><\tilde{Jf_l},\tilde{H}>$$
and $<J\nu,\tilde{H}>=0$ where J is the complex structure in
$C^{n+1}$, we have:
\begin{equation*}
\begin{split}
<Jf_k,H><B_{kl},H><Jf_l,H>=&<Je_\alpha,\tilde{H}><\tilde{B}_{\alpha
\beta},\tilde{H}><Je_\beta,\tilde{H}>\\
=&<Je_\alpha,\tilde{H}>^2<\tilde{B}_{\alpha \alpha},\tilde{H}>
\end{split}
\end{equation*}

Every term in it has explicit expression, thus:
$$\frac{1}{2}\frac{d^2}{dt^2}E_2(i_t)|_{t=0}=-4\int[d-(n+1)^2]^2+3[2(n+1)^3+(\sum\frac{1}{a_i^4})-3(n+1)d]d
v_g$$

Recall (11) and $t=2(n+1)(n+3)-d$, after tedious computation, we
get:
$$\frac{1}{2}\frac{d^2}{dt^2}E_2(i_t)|_{t=0}<0$$
i.e. the torus are unstable.

By (12), we easily get:
\begin{proposition}
Let M be a compact proper biharmonic submanifold of $CP^n$, if the
eigenvalues of its Weingarten operator are all nonnegative, then it
can not be stable.
\end{proposition}
Although it can't apply to our examples.

All the computations in above sections are valid in $QP^n$ and
$S^{4n+3}$, so the same method can generate lots of examples of
biharmonic submanifolds in $QP^n$ as well.

Finally, the author want to thank his advisor Prof. Dong for the
instructions, and Prof. Chao and Prof. Ji for valuable
conversations.

\newpage

\noindent Wei Zhang

\

\noindent School of Mathematical Sciences

\noindent Fudan University

\noindent Shanghai, 200433, P. R.China

\

\noindent Email address: 032018009@fudan.edu.cn


\begin{thebibliography}{99}

\bibitem[BB]{BB} The Bibliography of Biharmonic Maps, {http://beltrami.sc.unica.it/biharmoinc/}.

\bibitem[BMO]{BMO} A.Balmus, S.Montaldo and C.Oniciuc,
{Classification results for biharmonic submanifolds in spheres},
preprint.

\bibitem[CMO]{CMO} R.Caddeo, S.Montaldo, C.Oniciuc, {Biharmonic
submanifolds of $S^3$}. Internat. J. Math. {\bf 12} (2001), 867-876.

\bibitem[J1]{J1} G. Y. Jiang, {2-harmonic isometric immersions between Riemannian manifolds}, Chinese Ann. Math. Ser. A, {\bf 7} (1986), 130-144.

\bibitem[J2]{J2} G. Y. Jiang, {2-harmonic maps and their first and
second variational formulas}, Chinese Ann. Math. Ser. A, {\bf 7}
(1986), 389-402.

\bibitem[L] {L} B.Lawson, {Rigidity Theorem in Rank-1 symmetric
spaces}, J. Differential Geometry, {\bf 4} (1970), 349-357.

\bibitem[LO]{LO} E.Loubeau and C.Oniciuc, {On the biharmonic and
harmonic indices of the Hopf map}, Trans. Amer. Math. Soc., to
appear.

\bibitem[MO]{MO} S.Montaldo and C.Oniciuc, {A short survey on
biharmonic maps between Riemannian manifolds}, {\bf Proceedings of
the "II Workshop in Differential Geometry"}, Cordoba, June 2005, to
appear.

\bibitem[WL]{WL} Chuanxi Wu and Guanghan Li, {Geometry of
submanifolds(Chinese)}, Science Press, Beijing, 2002.



\end{thebibliography}
\end{document}